\documentclass[12pt]{article}
\usepackage{hyperref}
\usepackage{amsmath}
\usepackage{ntheorem}
\usepackage{mathtools}
\usepackage{epsf}
\usepackage{theorem}
\usepackage{indentfirst}
\usepackage{latexsym}
\usepackage{amssymb}
\usepackage[dvips]{graphicx}
\usepackage{flafter}
\usepackage{caption}
\usepackage{array}
\usepackage{textcomp}
\usepackage{supertabular}
\usepackage{longtable, booktabs}
\usepackage{hhline}
\usepackage{bigstrut,bigdelim}
\usepackage{rotating}
\usepackage{multicol}
\usepackage{multirow}
\usepackage{makeidx}
\usepackage{epsfig}
\usepackage{fancyhdr}
\newif\ifger
\gerfalse

\newtheorem{lemma}{Lemma}[section]
\newtheorem{theorem}{Theorem}
\newtheorem{definition}{Definition}[section]
\newtheorem{proposition}{Proposition}[section]
\newtheorem{corollary}{Corollary}[section]

\newtheorem{remark}{Remark}
\newtheorem{rem}{Remark on Table}

\begin{document}
\baselineskip=21pt

\title{ Block-Transitive Automorphism Groups of $2$-$(v,5,\lambda)$ Designs}

\author{Chuhan Lei, Xiaoqin Zhan\footnote{ Corresponding author. E-mail address:  zhanxiaoqinshuai@126.com}\\
\small \it  School of Science, East China JiaoTong University, \\
\small \it  Nanchang, 330013, P.R. China}

\date{}
\maketitle

\begin{abstract}
This paper investigates $2$-$(v,5,\lambda)$ designs $\mathcal{D}$ admitting a block-transitive automorphism group $G$. We first prove that if $G$ is point-imprimitive, then $v$ must be one of 16, 21, or 81. We further provide a complete classification of all such designs for $v=16$ and $v=21$.
Secondly, we demonstrate that if $G$ is point-primitive, then it must be of affine type, almost simple type, or product type. Additionally, we present a classification of pairs $(\mathcal{D},G)$ where $G$ is of product type.



\end{abstract}

\section{Introduction}

 \begin{definition}
  {\rm For positive integers $t< k< v-1$ and $\lambda$, a $t$-$(v,k,\lambda)$ design $\mathcal{D}$ is an incidence structure $(\mathcal{P},\mathcal{B})$ satisfying the following properties:

   \begin{itemize}

   \item
    $\mathcal{P}$ is a set of $v$ elements, called {\it points},
   \item
   $\mathcal{B}$ is a set of $b$ $k$-subsets of $\mathcal{P}$, called {\it blocks},
    \item
    Every $t$-subset of $\mathcal{P}$ is contained in exactly $\lambda$ blocks.
 \end{itemize}}
 \end{definition}
  Since all the blocks have the same size $k$, it follows that each point belongs to the same number $r$ of blocks.

An {\it automorphism} of $\mathcal{D}$ is a
permutation of the point set that preserves the block set. The set
of all automorphisms of $\mathcal{D}$ with the composition of permutations forms a group, denoted by $\rm {Aut}({\cal D})$. 
A subgroup $G$ of $\rm {Aut}({\cal D})$ is called
{\it point-primitive} if its  action  on the point set $\mathcal{P}$ is primitive, and {\it point-imprimitive} otherwise. A {\it flag} of $\mathcal{D}$ is a  point-block pair $(\alpha, B)$ with $\alpha \in B$, and $\cal D$ is called {\it flag-transitive} ({\it block-transitive}) if $G$ acts transitively on the set of flags (blocks).
A set of blocks of incidence structure $\mathcal{D}$ is called a set of {\it base blocks} with respect to an automorphism group $G$ of $\mathcal{D}$ if it contains exactly one block from each $G$-orbit on the block set. In particular, if $G$ is a block-transitive automorphism group of $\mathcal{D}$, then any block $B$ is a base block of $\mathcal{D}$.

If $G$ is a flag-transitive automorphism group of a $t$-design $\cal D$, then $G$ is clearly both block-transitive and point-transitive. In particular, as was proved by Ray-Chaudhuri and Wilson \cite{Ray} flag-transitivity implies point-primitivity in any $t$-$(v,k,\lambda)$ design where $t\ge3$. Some counterexamples given by Davies \cite{Dav1} show that this is no longer true for $t=2$.  It was proved by  Higman and McLaughlin \cite{HMcL} that if  $G$ is a flag-transitive automorphism group of a $2$-$(v,k,1)$ design then $G$ is point-primitive.  Cameron and Praeger \cite{Cameron} showed that for a flag-transitive, point-imprimitive $2$-$(v,k,\lambda)$ design, the number of points $v$ satisfies  $v\leq (k-2)^2$.
Montinaro in \cite{Montinaro} gave a classification of flag-transitive $2$-$(v,k,\lambda)$ designs with $v=k^2$ and $\lambda|k$. 

It follows from a result of  Delandtsheer and Doyen \cite{Delan} that if a $2$-$(v,k,\lambda) $ design is block-transitive and point-imprimitive, then $v\le(\binom{k}{2}-1)^2$ (see Lemma \ref{imn}).
It was shown in \cite{Cameron} that for a block-transitive and point-imprimitive $t$-design then $t\le3$. Also, they conjectured that $3$-$(v,k,\lambda)$ designs should also be further restricted: namely if such a 3-design is block-transitive and point-imprimitive, then $v\le \binom{k}{2}+1$ should hold. Mann and Tuan in \cite{Tuan1} proved this conjecture, and Tuan in \cite{Tuan2} constructed infinitely many examples of 3-designs attaining
this bound. It follows immediately that, for a block-transitive, point-imprimitive $2$-$(v,k,\lambda)$ design, the number $v$ of points cannot be too large relative to the block size $k$. Thus, if $\mathcal{D}$ is a $2$-$(v,k,\lambda)$ design with a block-transitive and point-imprimitive automorphism group $G$, then  $k\geq4$.

 Block-transitive 2-designs with small block size $k$ are of significant interest. For example, triple systems ($2$-$(v,3,\lambda)$ designs) have been extensively studied (see \cite{ZhanDing,ZhanPD}). For $k=4$, the second author, Zhou and Chen \cite{ZhanJCD} proved that a flag-transitive automorphism
group of a $2$-$(v,4,\lambda)$ design is point-primitive of affine, almost simple or product type,
and they classified all such designs with product type automorphism groups. The
almost simple case for $k=4$ or 5, when the socle is $PSL(2,q)$ or a sporadic
simple group, has been treated (see \cite{Ding5,ZhanDCC,Shen}). The 2-designs with block size 6 admitting
flag-transitive and point-imprimitive automorphism groups have been determined in \cite{ZhanZhou}.

Motivated by the work of the second author in \cite{ZhanDM}, this paper continues to investigate block-transitive
2-designs with small block size. In particular, we classify block-transitive and point-imprimitive 2-designs with block size 5. Our first result is as follows:
\begin{theorem}\label{th1}
Let $\mathcal{D}=(\mathcal{P},\mathcal{B})$ be a $2$-$(v,5,\lambda)$ design, admitting a block-transitive and point-imprimitive automorphism group $G$ that leaves a invariant nontrivial partition of $\mathcal{P}$ with $d$ classes of size $c$. Then the possible values for $(v,c,d)$ are:
$$(16,4,4),(21,3,7),(21,7,3),(81,9,9).$$ 
Moreover, 
\begin{enumerate}
  \item[\rm(i)]  If $(c,d)=(4,4)$ then $\lambda\in\{4_5,8_{15},12_2,16_{12},24_3,32_3,48_3,96_1,144_1\}$.
  \item[\rm(ii)]  If $(c,d)=(3,7)$ then $\lambda$ is one of the following: $$\{1_1,2_4,3_6,4_3,6_{26},8_8,12_{14},20_1,24_4,27_1,48_1,54_1,60_1,81_3,108_2,120_1,162_3,324_1,540_1\}.$$
  \item[\rm(iii)]  If $(c,d)=(7,3)$ then $\lambda\in\{7_3,14_5,21_6,42_5,49_1,98_1,147_1,196_1,245_1\}.$
  \item[\rm(iv)] If $(c,d)=(9,9)$ then $\lambda\le 40824.$
\end{enumerate}

\end{theorem}

\begin{remark}
\begin{enumerate}
{\rm
  \item[(i)] The notation $\lambda=a_n$  indicates that $\lambda=a$ and there are exactly  $n$ pairwise non-isomorphic block-transitive and point-imprimitive $2$-$(v,5,\lambda)$ designs. The base block and automorphism groups of each design are given in Tables \ref{BAC} and \ref{BA}.   
  \item[(ii)] All designs in this paper are constructed by software package {\sc Magma} \cite{Magma}. Since the transitive groups in {\sc Magma} are only given for degrees up to 47, thus we cannot yet determine all the block-transitive and point-imprimitive $2$-$(81,5,\lambda)$ designs.}
\end{enumerate}
\end{remark}

 Next, we consider the case where $G$ is block-transitive and point-primitive. Applying the O'Nan-Scott Theorem \cite{ONSC} for finite primitive permutation groups, we arrive at the following result:

\begin{theorem}\label{th2}
  Let $\mathcal{D}=(\mathcal{P},\mathcal{B})$ be a $2$-$(v,5,\lambda)$ design with a block-transitive and point-primitive automorphism group $G$. Then $G$ must be of affine type, almost simple type, or product type.
\end{theorem}

Existing research on block-transitive and point-primitive $2$-$(v,k,\lambda)$ designs has primarily focused on the affine and almost simple types. However, the product type has received significantly less attention. This paper focuses specifically on 2-$(v,5,\lambda)$ designs and investigates the case where the automorphism group is block-transitive, point-primitive, and of product type. Our analysis yields the following results:
\begin{theorem}\label{a1}
 Let $\mathcal{D}=(\mathcal{P},\mathcal{B})$ be a $2$-$(v,5,\lambda)$ design with a block-transitive point-primitive automorphism group $G$. Assume that $G$ has a product action on the set $\cal P$, then one of the following holds:
  \begin{enumerate}
    \item[\rm(i)] ${\rm Soc}(G)=A_9\times A_9$, and $\mathcal{D}$ is a $2$-$(81,5,\lambda)$ design with $\lambda \in\{5880,7056,14112\}$.
    \item[\rm(ii)] ${\rm Soc}(G)=PSL(2,8)\times PSL(2,8)$, and $\mathcal{D}$ is a $2$-$(81,5,\lambda)$ design with $$\lambda \in\{392,784,1176,1568,2352,4704,7056\}.$$
    \item[\rm(iii)] ${\rm Soc}(G)=A_{19}\times A_{19}$, and $\mathcal{D}$ is a $2$-$(361,5,3329280)$ design.
  \end{enumerate}
  
\end{theorem}


\section{Preliminaries}
We follow standard notation and terminology from design theory (see \cite{Handbook, Demb1968}) and group theory (see \cite{Wie1964}).
\begin{lemma}{\rm \cite[1.2, 1.9]{Handbook}} \label{L1}
The parameters $v, b, k, r$  and $\lambda$ of a $2$-design
satisfy the following conditions:
\begin{enumerate}
\item[\rm(i)] \, $vr=bk$.
\item[\rm(ii)] \, $\lambda(v-1)=r(k-1)$.
\end{enumerate}
\end{lemma}

From Lemma \ref{L1}(ii), it is easy to draw the following corollary: 
\begin{corollary}\label{a2}
Let $\mathcal{D}$ be a $2$-$(v,5,\lambda)$ design, then $r>\frac{\lambda v}{5} $.
\end{corollary}

Let $G_{\alpha}$ denote the stabilizer in $G$ of a point $\alpha\in\mathcal{P}$, and $G_B$ the setwise stabilizer in $G$ of a block $B\in\mathcal{B}$.

\begin{lemma}\label{divide}
  Let $\mathcal{D}$ be a $2$-$(v,5,\lambda)$ design with a block-transitive automorphism group $G$. Then $r$ divides $5\lambda d$, for all non-trivial subdegrees $d$ of $G$. Moreover, $v-1$ divides $20d$.
  \end{lemma}
 \textbf{Proof.}\quad
 Let $B$ be a block of $\mathcal{D}$ containing the point $\alpha$. Then
 $$[G:G_{\alpha B}]=[G:G_\alpha][G_\alpha:G_{\alpha B}]=[G:G_B][G_B:G_{\alpha B}].$$
 So $[G_\alpha:G_{\alpha B}]=\frac{r[G_B:G_{\alpha B}]}{5}$ by Lemma \ref{L1}(i). If $[G_B:G_{\alpha B}]=5$, then $G$ acts flag-transitively on $\cal D$ and so $r\mid \lambda d$ by \cite{Dav1}. If $[G_B:G_{\alpha B}]\leqslant 4$, then $\frac{r}{5}\mid [G_\alpha:G_{\alpha B}]$. 
 Let $\Gamma \neq \{\alpha\}$ be a non-trivial orbit of $G_\alpha$ on $\mathcal{P}$. Suppose that $G_\alpha$ has $t$ orbits on pencil $P(\alpha)$ (the set of blocks through $\alpha$), and denoted by $O_1,O_2,\dots,O_t$ where $1\leq t\leq 5$. Set $\mu _i=|\Gamma \cap B_i|$, where $B_i\in O_i$ $(1\leq i\leq t)$. Clearly, $0\leq \mu _i\leq 4$. Computing the size of  $\{(\beta ,B)|B\in P(\alpha),\beta \in \Gamma \cap B\}$ in two different ways, we get $$\sum_{i = 1}^{t} |O_i|\mu _i=\lambda |\Gamma|. $$
 So $\frac{r}{5}\mid \lambda|\Gamma | $ as $\frac{r}{5}\mid |O_i| $, thus $r\mid 5 \lambda |\Gamma|$. 
 
Then according to the Lemma \ref{L1}(ii), we have $\frac{\lambda(v-1)}{4} \mid 5\lambda |\Gamma|$. Because of this, the result is easily established. 
$\hfill\square$

In \cite{Delan}, Delandtsheer and Doyen proved that in a block-transitive point-imprimitive $2$-$(v,k,\lambda)$ design, $v$ is bounded by a function of $k$.
\begin{lemma}{\rm \cite[Theorem 5.1]{Delan}}\label{imn}
If $\mathcal{D}$ is a block-transitive point-imprimitive $2$-$(v,k,\lambda)$ design, then $v\le(\binom{k}{2}-1)^2$. 
\end{lemma}

The basis of our method is the following elementary result.
\begin{lemma}{\rm \cite[Proposition 1.1]{Delan}}\label{ele}
Let $\mathcal{D}=(\mathcal{P},\mathcal{B})$ be a $t-(v,k,\lambda)$ design, admitting a block-transitive group $G$. Let $H$ be a permutation group with $G\le H\le S_v$, and $\mathcal{B}^*=\mathcal{B}^H$ the set of images of blocks in $\cal B$ under $H$. Then $(\mathcal{P},\mathcal{B}^*)$ is a $t-(v,k,\lambda^*)$ design, for some $\lambda^*$, admitting the block-transitive automorphism group $H$.
\end{lemma}

\section{Proofs}
In this section, we focus on $2$-$(v,5,\lambda)$ designs with a block-transitive automorphism group $G$. Section 3.1 presents a complete classification of all pairs $(\mathcal{D},G)$, where $\cal D$ is a $2$-$(v,5,\lambda)$ design with $v\in\{16,21\}$ and $G$ is a block-transitive automorphism group acting point-imprimitive on the points of $\mathcal{D}$. Subsequently, in Section 3.2, we first apply a reduction to $G$ and then proceed to classify all pairs $(\mathcal{D},G)$ where $G$ acts point-primitive on the points of $\cal D$ and is of product type.

\subsection{Imprimitive Case}

Firstly, we are interested in the possibility that a block-transitive point-imprimitive group $G=S_c\wr S_d$ of an incidence structure $(\mathcal{P},\mathcal{B}^*)$. Suppose  the set of points is partitioned into $d$ non-trivial blocks of imprimitivity $\Delta_j$, $j=1,\ldots,d$, each of size $c$. So $v=cd$, with $c, d>1$.  Then the sizes of the intersections of each block with the imprimitivity classes
determine a partition of $k$, say $\mathbf{x} = (x_1, x_2, \ldots , x_d )$ with $x_1\ge x_2\ge \ldots \ge x_d $ and $\sum\limits_{i=1}^d x_i = k$. Denote this structure by $\mathfrak{D}(c,d,\mathbf{x})$. Set $b_t=\sum\limits_{i=1}^d x_i(x_i-1)\cdots(x_i-t+1)$. Note that $b_1=k$. By \cite[Proposition 2.2]{Cameron}, the following lemma holds.

\begin{lemma} \label{b2}
Let $\mathcal{D}^*=\mathfrak{D}(c,d,\mathbf{x})$. Then $\mathcal{D}^*$ is a $2$-design if and only if $$b_2=\sum_{i=1}^d x_i(x_i-1)=\frac{k(k-1)(c-1)}{(v-1)}.$$
\end{lemma}

In the following we always assume $v=cd$, with $c, d>1$.

\begin{proposition} \label{para}
 Let $\mathcal{D}^*=(\mathcal{P},\mathcal{B}^*)$ be a $2$-$(v,5,\lambda)$ design. If $G=S_c\wr S_d$ acts as  a block-transitive  point-imprimitive automorphism group of $\mathcal{D}^*$, then one of the following holds:
 \begin{enumerate}
  \item[\rm(i)] \, $(v,\lambda)=(16,144)$ and $c=d=4$.
  \item[\rm(ii)] \, $(v,\lambda)=(21,540)$ and $c=3, d=7$.
  \item[\rm(iii)] \, $(v,\lambda)=(21,245)$ and $c=7, d=3$.
  \item[\rm(iv)] \, $(v,\lambda)=(81,40824)$ and $c=d=9$.
  \end{enumerate}
\end{proposition}
\textbf{Proof.}\quad
 From Lemma \ref{imn}, we have $v\le81$.
Clearly, $b_2=\sum\limits_{i=1}^d x_i(x_i-1)$ is an even number with $x_i\le c$  for $1\le i\le d$.   By Lemma \ref{b2}, the 3-tuple $(v,c,d)$ is one of the following: 
$$(16,4,4),(21,3,7),(21,7,3),(81,9,9).$$

If $v=16$ with $c=d=4$. We have $\mathbf{x}=(x_1,x_2,x_3,x_4)=(2,2,1,0)$. Let $G=S_4\wr S_4$ acts imprimitively on points, with  $\Delta_1, \Delta_2, \Delta_3, \Delta_4$ representing the blocks of imprimitivity.  Consider a block  $B\in \mathcal{B}^*$ such that $|B\cap\Delta_i|=x_i$ for $1\le i\le4$. Since $G$ acts $4$-transitively on $\{\Delta_1, \Delta_2, \Delta_3, \Delta_4\}$, there exists a $g\in G$ such that $B^g=C$, where $C$ is a $5$-element subset of $\cal P$ with a partition $\mathbf{x}=(2,2,1,0)$. By block-transitivity of $G$, we have
 $$b=\binom{4}{2}\binom{4}{2}^2\binom{2}{1}\binom{4}{1}=1728.$$
  Applying Lemma \ref{L1}, we get $\lambda=144$. Moreover, there is a unique $2$-$(16,5,144)$ design admitting $G=S_4\wr S_4$ as its block-transitive  point-imprimitive automorphism group.

If $v=21$ with $c=3, d=7$. We have $\mathbf{x}=(x_1,x_2,x_3,x_4,x_5,x_6,x_7)=(2,1,1,1,0,0,0)$. Let $G=S_3\wr S_7$ acts imprimitively on points, with $\Delta_1, \Delta_2, \Delta_3, \Delta_4,\Delta_5, \Delta_6, \Delta_7$ representing  the blocks of imprimitivity. The block-transitivity of $G$ yields 
 $$b=\binom{7}{1}\binom{3}{2}\binom{6}{3}\binom{3}{1}^3=11340.$$
   By  Lemma \ref{L1}, we get $\lambda=540$. Furthermore, there exists a unique $2$-$(21,5,540)$ design admitting $G=S_3\wr S_7$ as its block-transitive  point-imprimitive automorphism group.
 
If $v=21$ with $c=7,d=3$. We have $\mathbf{x}=(x_1,x_2,x_3)=(3,1,1)$. Let $G=S_7\wr S_3$ acts imprimitively on points, with  $\Delta_1, \Delta_2, \Delta_3$ representing  the blocks of imprimitivity. By block-transitivity of $G$ we have
 $$b=\binom{3}{1}\binom{7}{3}\binom{7}{1}^2=5145.$$
  And $\lambda=245$ by  Lemma \ref{L1}. Moreover, a unique $2$-$(21,5,245)$ design exists admitting $G=S_7\wr S_3$ as its block-transitive  point-imprimitive automorphism group.
 
If $v=81$ with $c=d=9$. Then $\mathbf{x}=(x_1,x_2,x_3,x_4,x_5,\ldots,x_9)=(2,1,1,1,0,\ldots,0)$. Let $G=S_9\wr S_9$ acts imprimitively on points, with  $\Delta_1, \Delta_2, \Delta_3, \Delta_4,\Delta_5, \Delta_6, \Delta_7, \Delta_8, \Delta_9$  representing the blocks of imprimitivity. By block-transitivity of $G$ we have
 $$b=\binom{9}{1}\binom{9}{2}\binom{8}{3}\binom{9}{1}^3=13226976.$$
  And $\lambda=40824$ by  Lemma \ref{L1}. Moreover, there exists a unique $2$-$(21,5,40824)$ design admitting $G=S_9\wr S_9$ as its block-transitive  point-imprimitive automorphism group.
 $\hfill\square$

By Lemma \ref{ele}, we have the following result:
\begin{corollary}\label{value}
Let $\mathcal{D}=(\mathcal{P},\mathcal{B})$ be a  $2$-$(v,5,\lambda)$ design, admitting a block-transitive and point-imprimitive automorphism group $G$ that leaves a invariant non-trivial partition of $\mathcal{P}$ with $d$ classes of size $c$. Then $(v,c,d)$ is one of the following:
$$(16,4,4),(21,3,7),(21,7,3),(81,9,9).$$ 
\end{corollary}

In this study, we utilize the {\sc Magma} to determine all pairs $(\mathcal{D},G)$ where $\mathcal{D}$ is a $2$-$(v,5,\lambda)$ design with $v\in\{16,21\}$, and $G$ is a block-transitive point-imprimitive automorphism group. 

We begin by using the command ${\tt N:=TransitiveGroups(v)}$ to retrieve all transitive groups acting on $v$ points.  For $v=16$, there are 1954 transitive groups, of which 22 are primitive (see \cite[p829]{Handbook}). We  focus on the remaining 1932 imprimitive groups.
The command {\tt G:=N[i]} returns the $i$-th transitive group in the list of the {\sc Magma}-library of transitive groups with degree $16$. This command yields the transitive permutation representations of $G$ acting on the set $\mathcal{P}=\{1,2,3,\ldots,16\}$. For $v=21$, there are 164 transitive groups, of which 9 are primitive (see \cite[p830]{Handbook}). We focus our investigation on the remaining 155 imprimitive groups.

\begin{proposition}\label{pro1}
Let $\mathcal {D}=(\mathcal{P},\mathcal{B})$ be a $2$-$(v,5,\lambda)$ design with $v\in\{16,21\}$, and $G:=N[i]$ be a block-transitive point-imprimitive automorphism group of $\cal D$.
\begin{enumerate}
  \item[\rm(i)] If $v=16$ then $(\lambda,i)$ is given in Table \ref{TAB3}. 
  \item[\rm(ii)]  If $v=21$ then $(\lambda,i)$ is given in Table \ref{TAB2} for $(c,d)=(3,7)$, or $(\lambda,i)$ is given in Table \ref{TAB4} for $(c,d)=(7,3)$.
  \end{enumerate} 
\end{proposition}

\begin{table}[htbp]
 \centering
 \caption{ Parameters $(\lambda,i)$ for $v=16$}\label{TAB3}
 \begin{tabular}{c|llllllllllllll}
 \toprule
 $\lambda$&&&&&&&$i$&&&&&&&\\
 \midrule
 4&63&64&183&184&185&195&430&433&&&&&&\\
 \hline
 8&183&184&185&194&195&419&420&430&431&433&436&\\
 &440&731&745&746&769&776&1063&1073&\\
 \hline
 12&414&709&710&1033&&&&&&&&&&\\
 \hline
 16&419&420&430&431&433&436&440&731&745&746&749&\\
 &776&1042&1059&1063&1070&1073&1074&1314&1317&1545&1690&&\\
 \hline
 24&709&710&1029&1297&&&&&&&&&\\
 \hline
 32&745&746&776&1063&1070&1073&1074&1314&1317&1690&&&&\\
 \hline
 48&1029&1033&1297&1492&1649&1650&1652&1749&1750&1789&1790&\\
 &1797&1834&\\
 \hline
 96&1649&1650&1750&1789&1790&1834&&&&&&&&\\
 \hline
 144&1859&1876&1877&1881&1895&1899&1900&1901&1904&1914&1917&\\
 &1918&1919&1920&1929&1930&1937&1941&1942&1943&1947&\\
 \bottomrule
 \end{tabular}
\end{table}

\begin{table}[htbp]
 \centering
 \caption{ Parameters $(\lambda,i)$ for $(c,d)=(3,7)$ }\label{TAB2}
 \begin{tabular}{c|llllllllllllll}
 \toprule
 $\lambda$&&&&&&&$i$&&&&&&&\\
 \midrule
 1&1 &2&7&11&&&&&&&&&&\\
 \hline
 2&3&4&5&8&9&10&15&&\\
 \hline
 3&7 &11&&&&&&&&&&&&\\
 \hline
 4& 14&22&27&&&&&&&&&&&\\
 \hline
 6&9&10&11&15&&&&&&&&&\\
 \hline
 8& 14&&&&&&&&&&&&&\\
 \hline
 12&15&22&27&&&&&&&&&&&\\
 \hline
 20&44&56&57&58&74&&&&&&&&&\\
 \hline
 24&22&27&&&&&&&&&&&\\
 \hline
 27&39&50&59&60&61&79&80&86&100&105&112&114&\\
 &117&120&123&127&137&\\
 \hline
 48&27&&&&&&&&&&&\\
 \hline
 54&51&52&68&75&76&77&78&81&92&93&98&99&&\\
 &107&122&124&131&133&134&142&&&& \\
 \hline
 60&44&56&&&&&&&&&&&&\\
 \hline
 81&59 &60&79&86&100&114&117&127&137&&&&&\\
 \hline
 108&104 &111&113&115&118&119&125&126&135&136&145&147&&\\
 \hline
 120&57&58&74&&&&&&&&&\\
 \hline
 162&77&78&81&92&98&99&107&133&134&142&\\
 \hline
 324&104&111&113&115&118&119&125&126&135&136&145&147&\\
 \hline
 540&121&128&129&130&132&138&139&140&141&144&148&149&\\
 &150&151&152&&\\
 \bottomrule
 \end{tabular}
\end{table}

\begin{table}[htbp]
 \centering
 \caption{ Parameters $(\lambda,i)$ for $(c,d)=(7,3)$}\label{TAB4}
 \begin{tabular}{c|llllllllllllll}
 \toprule
 $\lambda$&&&&&&&$i$&&&&&&&\\
 \midrule
 7&12&13&18&21&26&&&&&&&\\
 \hline
 14&16&17&19&23&24&25&29&&\\
 \hline
 21&21&24&25&26&29&&&&&&&&&\\
 \hline
 42&24&25&29&&\\
 \hline
 49&28&32&34&35&40&48&49&65&73&90&143&146&&\\
 \hline
 98&30&31&36&37&41&42&43&45&46&47&53&54&\\
 &55&62&63&64&66&69&70&71&72&82&83&84&\\
 &87&88&89&94&95&96&97&101&102&106&108&109&\\
 &110&116&\\
 \hline
 147&34&35&40&41&42&43&48&49&53&54&55&63&\\
 &64&65&66&69&70&71&72&73&82&83&84&87&\\
 &88&89&90&94&95&96&97&101&102&106&108&109&\\
 &110&116&&\\
 \hline
 196&143&146&&&&&&&&\\
 \hline
 245&153&154&155&156&157&158&159&160&161&162&&\\
 \bottomrule
 \end{tabular}
\end{table}
{\bf Proof:} We present a proof for the case $(v,c,d) = (16, 4, 4)$. The remaining cases can be established using similar arguments.
  Let $H=S_4\wr S_4$, then there exists  a unique $2$-$(16,5,144)$ design $\mathcal{D}^*=(\mathcal{P},B_1^H)$ by Lemmas \ref{para}, where $B_1\subsetneq \mathcal{P}$ and  the sizes of the intersections of $B_1$ with the imprimitivity classes is $(2,2,1,0)$. Let $B$ be a base block of $\cal D$. Since the group $S_4\wr S_4$ is a maximal subgroup of $S_{16}$, any automorphism group $G$ of $\cal D$ must satisfy $G\le H$. Therefore, by iterating through all elements in $B_1^H$ and all potential automorphism groups $G:=N[i]$ from {\sc Magma}'s library of transitive groups with degree 16,  we can efficiently determine the existence of the design using the command {\tt Design<2,16|$B^G$>}.    $\hfill\square$

To verify the non-isomorphism of the designs constructed in Proposition \ref{pro1}, we utilize the {\sc Magma}-command ${\tt IsIsomorphic(D1,D2)}$. This comparison yields the following results: 
  \begin{proposition}\label{pro2} Up to isomorphism:  
    \begin{enumerate}
      \item[\rm(i)] \,There are $45$ different $2$-$(16,5,\lambda)$ designs $\cal D$ which admitting a block-transitive point-imprimitive automorphism group $G:=N[i]$, the pairs $(\mathcal{D},G)$ listed in Table \ref{BAC}. 
      \item[\rm(ii)] \,There are $106$ different $2$-$(21,5,\lambda)$ designs $\cal D$ which admitting a block-transitive point-imprimitive automorphism group $G:=N[i]$, the pairs $(\mathcal{D},G)$ listed in Table \ref{BA}.
      \end{enumerate}
  \end{proposition}

 \begin{longtable}{cllc|cllc}
  \caption{ Base block and automorphism group of $2$-$(16,5,\lambda)$}\label{BAC}\\ \hline
  \endfirsthead
  \multicolumn{8}{l}{(Continued)}\\
  \hline
  \endhead
  \hline
  \multicolumn{8}{r}{(Continued on next page)}\\
  \endfoot\hline
  \endlastfoot
  $\lambda$&$Nr$ & base block $B$ & $G$ & $\lambda$&$Nr$ & base block $B$ & $G$\\
  \hline
  4&5&$ \{2,3,12,13,14\}$&$N[63]$&&&$\{3,7,8,10,11\}$&$N[419]$\\
    &&$ \{5,7,9,12,14\}$&$N[63]$&&&$ \{2,4,10,11,14\}$&$N[430]$\\
    &&$ \{1,3,10,11,15\}$&$N[63]$&&&$\{2,4,7,13,15\}$&$N[430]$\\
    &&$ \{1,2,5,8,11\}$&$N[63]$&&&$\{7,8,10,11,14\}$&$N[430]$\\
    &&$ \{6,7,9,10,14\}$&$N[64]$&&&$\{5,7,9,12,14\}$&$N[430]$\\
  
  8&15&$ \{7,8,10,11,14\}$&$N[183]$&&&$\{2,4,10,11,14\}$&$N[433]$\\
    &&$  \{2,4,6,13,14 \}$&$N[183]$&&&$\{2,4,7,13,15\}$&$N[433]$\\
    &&$  \{5,7,12,14,15\}$&$N[183]$&&&$\{7,8,10,11,14\}$&$N[433]$\\
    &&$  \{1,4,6,14,15 \}$&$N[184]$&&&$\{5,7,9,12,14\}$&$N[433]$\\
    &&$  \{2,4,10,11,14\}$&$N[184]$&&&$\{3,7,8,10,11\}$&$N[440]$\\
    &&$  \{2,4,7,13,15 \}$&$N[184]$&&&$\{7,8,11,14,16\}$&$N[776]$\\
    &&$  \{4,6,8,15,16\}$&$N[184]$&24&3&$\{7,8,10,11,14 \}$&$N[709]$\\
    &&$  \{4,7,8,9,10 \}$&$N[184]$&&&$ \{7,8,10,11,14 \}$&$N[710]$\\
    &&$  \{7,8,10,11,14\}$&$N[184]$&&&$\{6,7,9,10,14\}$&$N[1029]$\\
    &&$  \{5,7,9,12,14 \}$&$N[184]$&32&3&$\{1,4,6,14,15\}$&$N[745]$\\
    &&$  \{1,3,6,13,16\}$&$N[184]$&&&$\{5,7,9,12,14 \}$&$N[745]$\\
    &&$  \{6,7,11,12,16\}$&$N[184]$&&&$\{2,4,6,13,14 \}$&$N[776]$\\
    &&$  \{5,7,9,12,14 \}$&$N[195]$&48&3&$\{1,4,6,14,15\}$&$N[1029]$\\
    &&$  \{6,8,9,12,14\}$&$N[195]$&&&$\{6,8,9,13,15\}$&$N[1029]$\\
    &&$  \{4,5,8,13,14 \}$&$N[769]$&&&$\{7,8,10,11,14\}$&$N[1033]$\\

  12&2&$  \{7,8,10,11,14\}$&$N[414]$&96&1&$ \{8,9,12,13,15\}$&$N[1649]$\\
    &&$  \{1,3,5,7,12\}$&$N[414]$&144&1&$\{1,4,6,14,15\}$&$N[1859]$\\

  16&12&$\{1,4,6,14,15\}$&$N[419]$&&&&\\
   \end{longtable}

   \begin{longtable}{cllc|cllc}
    \caption{ Base block and automorphism groups of $2$-$(21,5,\lambda)$ design}\label{BA}\\ 
    \hline
    \endfirsthead
    \multicolumn{8}{l}{(Continued)}\\
    \hline
    \endhead
    \hline
    \multicolumn{8}{r}{(Continued on next page)}\\
    \endfoot\hline
    \endlastfoot
    $\lambda$&$Nr$ & base block $B$ & $G$ & $\lambda$&$Nr$ & base block $B$ & $G$\\
    \hline
    1&1&$ \{4,5,8,12,16\}$&$N[1]$&&&$\{1,6,9,11,12 \}$&$N[15]$\\
    
    2&4&$ \{4,6,7,10,18\}$&$N[3]$&&&$\{6,7,10,16,17\}$&$N[15]$\\
      &&$ \{3,4,7,9,10 \}$&$N[4]$&&&$\{4,6,8,10,13 \}$&$N[15]$\\
      &&$ \{5,12,16,20,21\}$&$N[4]$&&&$\{1,2,4,7,18 \}$&$N[15]$\\
      &&$ \{1,11,15,17,18\}$&$N[4]$&&&$\{4,9,12,20,21\}$&$N[15]$\\
    
    3&6&$ \{2,9,11,16,17\}$&$N[7]$&&&$\{5,10,12,18,21\}$&$N[15]$\\
     &&$  \{1,7,8,14,18 \}$&$N[7]$&&&$\{1,9,13,16,17 \}$&$N[15]$\\
     &&$  \{7,12,14,15,16\}$&$N[7]$&&&$\{5,12,14,16,17 \}$&$N[15]$\\
     &&$  \{1,3,7,16,21 \}$&$N[7]$&&&$\{1,12,13,14,17\}$&$N[22]$\\
     &&$  \{7,8,11,14,16 \}$&$N[7]$&&&$\{1,3,6,11,19\}$&$N[22]$\\
     &&$  \{1,13,15,18,20 \}$&$N[7]$&&&$\{2,7,11,12,14\}$&$N[22]$\\
    
    4&3&$\{1,6,10,13,15\}$&$N[14]$&&&$\{6,10,12,17,21\}$&$N[22]$\\
      &&$\{4,6,9,10,18 \}$&$N[22]$&14&5&$\{3,10,11,13,18\}$&$N[16]$\\
      &&$\{7,14,16,19,21\}$&$N[22]$&&&$\{3,5,6,14,18\}$&$N[16]$\\
    
    6&26&$ \{1,2,7,10,15 \}$&$N[9]$&&&$\{1,2,6,11,20\}$&$N[16]$\\
     &&$\{4,7,9,11,14\}$&$N[9]$&&&$\{3,10,11,13,18\}$&$N[17]$\\
     &&$\{9,10,11,17,19\}$&$N[9]$&&&$\{6,8,9,13,15\}$&$N[17]$\\
     &&$\{4,10,18,19,21 \}$&$N[9]$&20&1&$\{2,4,10,12,14\}$&$N[44]$\\
     &&$\{3,8,11,19,20\}$&$N[9]$&21&6&$\{1,4,6,14,15\}$&$N[21]$\\
     &&$\{2,3,7,11,21 \}$&$N[9]$&&&$\{6,8,9,13,15\}$&$N[21]$\\
     &&$\{1,8,9,12,13\}$&$N[9]$&&&$\{1,2,7,10,15\}$&$N[21]$\\
     &&$\{9,10,12,13,18\}$&$N[9]$&&&$\{6,13,15,16,19\}$&$N[21]$\\
     &&$\{1,3,9,15,18\}$&$N[9]$&&&$\{3,12,13,14,18\}$&$N[21]$\\
     &&$\{7,8,14,18,21 \}$&$N[9]$&&&$\{4,9,13,14,21\}$&$N[21]$\\
     &&$\{2,7,11,13,14 \}$&$N[9]$&24&4&$\{5,6,8,11,14\}$&$N[22]$\\
     &&$\{2,7,8,10,15\}$&$N[9]$&&&$\{6,7,10,16,17\}$&$N[27]$\\
     &&$\{2,3,8,12,21 \}$&$N[9]$&&&$\{2,7,11,12,14\}$&$N[27]$\\
     &&$\{2,9,11,19,21\}$&$N[10]$&&&$\{6,10,12,17,21\}$&$N[27]$\\
     &&$\{1,3,4,18,20 \}$&$N[10]$&27&1&$\{2,7,11,13,15\}$&$N[39]$\\
     &&$\{7,10,12,18,19\}$&$N[10]$&42&5&$\{1,4,6,14,15\}$&$N[24]$\\
     &&$\{3,8,11,19,20 \}$&$N[11]$&&&$\{6,8,9,13,15\}$&$N[24]$\\
     &&$\{1,3,6,9,20\}$&$N[11]$&&&$\{1,3,6,9,20\}$&$N[24]$\\
     &&$\{1,6,9,11,12\}$&$N[11]$&&&$\{2,3,7,11,21\}$&$N[24]$\\
     &&$\{6,7,10,16,17 \}$&$N[11]$&&&$\{1,6,7,12,19\}$&$N[24]$\\
     &&$\{4,6,8,10,13\}$&$N[11]$&48&1&$\{2,9,11,16,17\}$&$N[27]$\\
     &&$\{1,2,4,7,18 \}$&$N[11]$&49&1&$\{3,10,11,13,18\}$&$N[28]$\\
     &&$\{4,9,12,20,21\}$&$N[11]$&54&1&$\{2,7,11,13,15\}$&$N[51]$\\
     &&$\{5,10,12,18,21\}$&$N[11]$&60&1&$\{4,6,12,16,19\}$&$N[44]$\\
     &&$\{1,9,13,16,17 \}$&$N[11]$&81&3&$\{4,6,12,16,19\}$&$N[59]$\\
     &&$\{5,12,14,16,17 \}$&$N[11]$&&&$\{2,4,6,8,14\}$&$N[59]$\\
    
    7&3&$ \{3,10,11,13,18\}$&$N[12]$&&&$\{10,13,18,19,20\}$&$N[59]$\\
       &&$ \{3,10,11,13,18\}$&$N[13]$&98&1&$\{3,10,11,13,18\}$&$N[30]$\\
       &&$ \{6,8,9,13,15\}$&$N[13]$&108&2&$\{6,12,13,15,17\}$&$N[104]$\\

    8&8&$\{1,6,9,11,12 \}$&$N[14]$&&&$\{6,13,15,16,19\}$&$N[104]$\\
      &&$\{8,12,14,19,20\}$&$N[14]$&120&1&$\{3,8,11,19,20\}$&$N[57]$\\
      &&$\{3,5,6,14,18 \}$&$N[14]$&147&1&$\{1,4,6,14,15\}$&$N[34]$\\
      &&$\{6,9,15,16,18 \}$&$N[14]$&162&3&$\{4,6,12,16,19\}$&$N[77]$\\
      &&$\{10,13,16,17,21\}$&$N[14]$&&&$\{2,4,6,8,14\}$&$N[77]$\\
      &&$\{6,11,14,20,21\}$&$N[14]$&&&$\{10,13,18,19,20\}$&$N[77]$\\
      &&$\{1,3,11,15,18\}$&$N[14]$&196&1&$\{1,4,6,14,15\}$&$N[143]$\\
      &&$\{1,2,6,7,18\}$&$N[14]$&245&1&$\{1,4,6,14,15\}$&$N[153]$\\
    12&14&$\{3,8,11,19,20\}$&$N[15]$&324&1&$\{4,6,12,16,19\}$&$N[104]$\\
      &&$\{1,3,6,9,20\}$&$N[15]$&540&1&$\{4,6,12,16,19\}$&$N[121]$\\
  \end{longtable}

\begin{rem}
{\rm  The number $n$ listed in column ``Nr"  means that there are  $n$ pairwise non-isomorphic $2$-$(v,5,\lambda)$ designs admitting a block-transitive point-imprimitive automorphism group $G:=N[i]$.}
\end{rem}

Corollary \ref{value} and Propositions \ref{pro1}-\ref{pro2} complete the proof of Theorem \ref{th1}.  $\hfill\square$

\subsection{Primitive Case}

In this section, we give the proof of {Theorem} \ref{th2}.

\begin{proposition}\label{pro3} Let $\mathcal{D}=(\mathcal{P},\mathcal{B})$ be a $2$-$(v,5,\lambda)$ design and $G$ be a block-transitive point-primitive automorphism group of $\mathcal{D}$. Then $G$ is not of twisted wreath product type.
\end{proposition}
  {\bf Proof.}\, Suppose that $G$ is of twisted wreath product type. Let $N=T_1\times\cdots\times T_m $, where each $T_i\cong T$  $(i\in\{1,2,\ldots,m\})$ is a nonabelian simple group and $m\geq2$. Let $\alpha\in\mathcal{P}$. Then $G=N\rtimes G_{\alpha}$. $N$ is regular on $\mathcal{P}$ and $v=|T|^{m}$. Define subsets $\Gamma_i=\alpha^{T_i}$ for all $i\in\{1,2,\ldots,m\}$. Clearly, 
$|\Gamma_i|=|T_i:(T_i)_\alpha|=|T|$ and $\Gamma_i\cap\Gamma_j=\{\alpha\}$. Moreover, we have $$|\cup_{i=1}^m \Gamma_i-\{\alpha\}|=m(|T|-1).$$ 
Since $N\trianglelefteq G$, $G$ permutes $\{T_1,T_2,\ldots,T_m\}$ and so does $G_\alpha$. Thus, $\cup_{i=1}^m\Gamma_i \setminus\{\alpha\}$ is invariant under the action of $G_\alpha$. This implies that $\cup_{i=1}^m\Gamma_i \setminus\{\alpha\}$ is the union of some non-trivial $G_\alpha$-orbits. Consequently, $r$ divides $5\lambda m(|T|-1)$. From Corollary \ref{a2}, we get $$\frac{\lambda v}{5}<r<5\lambda m(|T|-1).$$ Hence, $|T|^m<25m(|T|-1)$. As $T$ is a nonabelian simple group, we have $$60^m\leq |T|^m<25m(|T|-1).$$ It follows that $m=1$, contradicting the assumption that $m\geq 2$.
  $\hfill\square$

\begin{proposition}\label{pro4} Let $\mathcal{D}=(\mathcal{P},\mathcal{B})$ be a $2$-$(v,5,\lambda)$ design and  $G$ be a block-transitive automorphism group of $\mathcal{D}$. Then $G$ is not of simple diagonal type.
\end{proposition}
  {\bf Proof.}\, Suppose that $G$ is of simple diagonal type. Let $N=T_1\times\cdots\times T_m $, where each $T_i\cong T$ $(i\in\{1,2,\ldots,m\})$ is a nonabelian simple group and $m\geq2$. Then $N\unlhd G\leq W$, where $$W=\{(a_1,a_2,\ldots,a_m)\pi \mid a_i\in {\rm Aut}(T), \pi \in S_m, a_i\equiv a_j (\bmod \,{\rm Inn}(T))\}.$$ For a point $\alpha\in\mathcal{P}$, we have $$T\cong \{(\alpha,\ldots,\alpha)\mid\alpha\in T\}=N_{\alpha}\leq G_{\alpha}\leq W_{\alpha},$$ where $$W_{\alpha}=\{(\alpha,\ldots,\alpha)\pi\mid\alpha\in {\rm Aut}(T),\pi \in S_m\}\cong {\rm Aut}(T)\times S_m.$$
  Clearly, $(T_i)_{\alpha}=T_i\cap N_{\alpha}=1$, $v=|N:N_{\alpha}|=|T|^{m-1}$. Define $\Gamma_i={\alpha}^T_i$ for $i\in\{1,2,\ldots,m\}$. We have $\Gamma_i\cap \Gamma_j=\{\alpha\}$ for distinct $i$, $j$. Choose an orbit $X$ of $G_{\alpha}$ in $\mathcal{P}\setminus\{\alpha\}$ such that $|X\cap \Gamma_1|=\mu _1\neq 0$. Let $m_1=|G_{\alpha}:N_{G_{\alpha}}(T_1)| $ and $ \mathcal{T}={T_1,T_2,\ldots,T_m}$. Since $G_{\alpha}\lesssim {\rm Aut}(T)\times S_m$ and $G^\mathcal{T}$ is transitive on $\mathcal{T}$, $m_1\leqslant m$. Thus $$d=|X|=m_1\mu _1\leq m|T|.$$
  As $r$ divides $5\lambda|X|$ and so $\frac{\lambda v}{5}<r<5\lambda m|T|$. Thus $|T|^{m-1}<25m|T|$. Since $T$ is a nonabelian simple group, we obtain that $$60^{m-2}\leq |T^{m-2}|<25m.$$ It follows that $m\leq 3$.
  Assume that $m=2$ or $m=3$. Hence, from the proof of {Lemma} \ref{divide}, we obtain that $r$ divides $5|G_\alpha|$, so we have $r$ divides $5|T||{\rm Out}(T)|m!$. From {Lemma} \ref{L1}, we obtain $r=\frac{\lambda}{4}(|T|^{m-1}-1)$. As $(|T|,|T|^{m-1}-1)=1$, which implies $r$ divides $\frac{5\lambda}{4}|{\rm Out}(T)|m!$. By {Corollary} \ref{a2}, we have $\lambda|T|^{m-1}<5r<\frac{25\lambda}{4}|{\rm Out}(T)|m!$, therefore we get $$4|T|^{m-1}<25|{\rm Out}(T)|m!.$$ 
  
  It follows that $|T|<\frac{25}{2}|{\rm Out}(T)| $ when $m=2$, or $|T|^2<\frac{75}{2}|{\rm Out}(T)| $ when $m=3$. 
  However, similar to the proof of \cite[Lemma 2.3]{Zhu}, it easily shows that there exists no group $T$ fulfilling these two inequalities. $\hfill\square$

\begin{proposition}\label{pro5} Let $\mathcal{D}=(\mathcal{P},\mathcal{B})$ be a $2$-$(v,5,\lambda)$ design,  and $G$ be a block-transitive automorphism group of $\mathcal{D}$. If $G$ is point-primitive of product type, then one of the following holds:
  \begin{enumerate}
    \item[\rm(i)] ${\rm Soc}(G)=A_9\times A_9$, $v=81$  and $\lambda \in\{5880,7056,14112\}$.
    \item[\rm(ii)]  ${\rm Soc}(G)=PSL(2,8)\times PSL(2,8)$,  $v=81$ and $$\lambda \in\{392,784,1176,1568,2352,4704,7056\}.$$
    \item[\rm(iii)] ${\rm Soc}(G)=A_{19}\times A_{19}$,  $v=361$ and $\lambda=3329280$.
  \end{enumerate}
\end{proposition}
{\bf Proof.}\, Suppose that $G$ has a product action on $\mathcal{P}$. Then there is a group $K$ with a primitive action (of almost simple or
diagonal type) on a set $\Gamma$ of size $v_{0}\geq5$, where
\begin{equation*}
  \mathcal{P}=\Gamma^{m}, G\leq K^{m}\rtimes S_{m}=K\wr  S_{m} \ \text{and}\ \ m\geq2.
\end{equation*}

 Let $H=K\wr  S_{m}$, and let $S_{m}$ act on $M=\{1,2,\ldots,m\}$. Let $\alpha$ and $\beta$ be two arbitrary points of $\cal P$. Clearly, $|\beta^{G_{\alpha}}|=[G_{\alpha}:G_{\alpha\beta}]$ is a non-trivial subdegree of $G$. By Lemma \ref{divide}, we have
 $r$ divides $5\lambda[G_{\alpha}:G_{\alpha\beta}]$, thus $[G_{\alpha}:G_{\alpha\beta}]\geq\frac{v-1}{20}$. As $G\le H=K\wr  S_{m}$, it follows that $G_{\alpha}\le H_{\alpha}$ and $|\beta^{G_{\alpha}}|\le|\beta^{H_{\alpha}}|$. Thus,
\begin{equation}\label{gongshi1} [H_{\alpha}:H_{\alpha\beta}]\geq[G_{\alpha}:G_{\alpha\beta}]
\geq\frac{v-1}{20}.
\end{equation}
Let $\alpha=(\gamma,\gamma,\ldots,\gamma)$, $\gamma\in\Gamma$,
$\beta=(\delta,\gamma,\ldots,\gamma)$, $\gamma\neq\delta\in\Gamma$
and $B\cong K^{m}$ be the base group of $H$. Then
$B_{\alpha}=K_{\gamma}\times\cdots\times K_{\gamma}$,
$B_{\alpha\beta}=K_{\gamma\delta}\times K_{\gamma}\times
\cdots\times K_{\gamma}$. Now $H_{\alpha}=K_{\gamma}\wr S_{m}$,
and $H_{\alpha\beta}=K_{\gamma\delta}\times(K_{\gamma}\wr S_{m-1})$. Suppose that $K$ has rank $s$ (the number of orbits of a point stabilizer in $K$) on $\Gamma$  with $s\geq2$.  We can choose a  $\delta\in \Gamma$ satisfying $[K_{\gamma}:K_{\gamma\delta}]\leq\frac{v_{0}-1}{s-1}$, so that
\begin{equation*}
[H_{\alpha}:H_{\alpha\beta}]=\frac{|H_{\alpha}|}
{|H_{\alpha\beta}|}
=\frac{|K_{\gamma}|^{m}\cdot m!}{|K_{\gamma\delta}||K_{\gamma}|^{m-1}\cdot(m-1)!}\leq m\frac{v_{0}-1}{s-1}
\end{equation*}
and hence by Equation (\ref{gongshi1}),
\begin{equation*}
  \frac{v-1}{20}\leq[G_{\alpha}:
  G_{\alpha\beta}]\leq m\frac{v_{0}-1}{s-1}.
\end{equation*}
So
\begin{equation}\label{gongshi2}
  \frac{v_{0}^{m}-1}{v_{0}-1}\leq \frac{20m}{s-1}\leq20m.
\end{equation}
We infer that $m=2$ with $s\in\{2,3,4,5,6,7\}$ or $m=3$ with $s=2$.

First, consider the case $m=3$ with $s=2$. Equation (\ref{gongshi2}) implies $5\leq v_0\leq 7$. The group $H=K\wr S_3$ has rank 4 with subdegrees 1, $3(v_{0}-1)$, $3(v_{0}-1)^2$ and $(v_{0}-1)^{3}$ on  $\mathcal{P}=\Gamma\times\Gamma\times\Gamma$. Since $G\leq H$, each subdegree of $H$ must be the sum of some subdegrees of $G$,  by Lemmas \ref{L1} and \ref{divide} we have $v_0^3-1$ dividing $60(v_{0}-1)$, which is a contradiction.

 Now consider the case $m=2$.
 
   If $4\leq s\leq 7$, then Equation (\ref{gongshi2}) determines the values of $v_0$, which are listed in Table \ref{TAB123456}. No group $K$ with a primitive action (of almost simple or diagonal type) exists on $\Gamma$ for any of these $v_0$ values.
  \begin{table}[htbp]
    \centering
    \caption{The values of  $v_0$  for different $s$}\label{TAB123456}
    \begin{tabular}{c|llllllllllllll}
    \toprule
    $s$&&&&&&&$v_0$&&&&&&&\\
    \midrule
    4&5&6&7&8&9&10&11&12&&&&&&\\
    \hline
    5&5&6&7&8&9&&\\
    \hline
    6&5&6&7&&\\
    \hline
    7&5&\\
    \bottomrule
    \end{tabular}
   \end{table}

  If $s=3$ then $5\le v_0\le19$ by Equation (\ref{gongshi2}). We get $v_0=10$ or $15$ as $K$ is primitive group with rank $3$ on $\Gamma$.
Moreover, if $v_0=10$, $H=K\wr S_2$  has rank 6 with subdegrees 1, 6, 9, 12, 36, 36 on $\mathcal{P}=\Gamma\times\Gamma$, which contradicts to Lemma \ref{divide}.   
  If $v_0=15$, $H=K\wr S_2$  has rank 6 with subdegrees 1, 12, 16, 36, 64, 96 on $\mathcal{P}=\Gamma\times\Gamma$, which contradicts to Lemma \ref{divide}. 
  
  If $s=2$ then $5\le v_0\le39$ by Equation (\ref{gongshi2}). $K$ acts 2-transitively on $\Gamma$. Therefore, $H=K\wr S_{2}$ has rank 3 with subdegrees 1, $2(v_{0}-1)$, $(v_{0}-1)^{2}$ on  $\mathcal{P}=\Gamma\times\Gamma$. Note that $G\leq H$,  each subdegree of $H$ must be the sum of some subdegrees of $G$.  By Lemma \ref{divide}, $v_0^2-1$ divides $40(v_0-1)$, which implies that $v_0+1$ divides 40,  leading to $v_0\in\{7,9,19,39\}$.

Suppose that $C$ is a 5-elements of $\mathcal{P}=\Gamma\times\Gamma$. Set  $\mathcal{B}=C^G$, and
\begin{center}
$\varPhi =\{B\in \mathcal{B}\mid (1,1), (1,2)\in B\}$,
$\varPsi =\{B\in\mathcal{B}\mid (1,1), (2,2)\in B\}.$
\end{center}
Then $\mathcal{D}=(\mathcal{P},\mathcal{B})$ is a $2$-$(v_0^2,5,\lambda)$ design if and only if $|\varPhi|=|\varPsi|=\lambda$ as $G$ has rank 3 on $\cal P$. Assume that ${\rm Soc}(G)=T\times T$. There are two possibilities here.

{\bf (a) $T$ is 5-transitive group of $\Gamma$}

  Let $C=\{(1,1),(1,2),a,b,c\}$, we consider all the possibilities of $\varPhi$ and $\varPsi$ enumerated in the Table \ref{BACaaa} below.
                                                                                                     
  \begin{longtable}{cllll}
  \caption{All the possibilities of $|\varPhi|$ and $|\varPsi|$}\label{BACaaa}\\ \hline
  \endfirsthead
  \hline
  \endhead
  \hline
  \endfoot\hline
  \endlastfoot
   Case&$\{a,b,c\}$  & $|\varPhi |$  & $|\varPsi |$  &\\
  \hline
  1&$\{(1,3),(1,5),(1,7)\}$ & $\frac{1}{6}(v_0-2)(v_0-3)(v_0-4) $ & $0$ &\\
  \hline
  2&$\{(1,7),(5,1),(5,2)\}$ & $6(v_0-1)(v_0-2)$ & $8(v_0-2)$ &\\
  \hline
  3& $\{(3,2),(4,2),(5,2)\}$ & $\frac{7}{3}(v_0-1)(v_0-2)(v_0-3) $ & $2(v_0-2)(v_0-3)$ &\\
  \hline
  4 &$\{(5,3),(5,5),(5,6)\}$ & $\frac{2}{3}(v_0-1)(v_0-2)(v_0-3)(v_0-4) $ & $2(v_0-2)(v_0-3)(v_0-4)$ &\\
  \hline
  5 & $\{(1,3),(3,3),(6,3)\}$ & $\frac{3}{2}(v_0-1)(v_0-2)^2$ & $2(v_0-2)^2$ &\\
  \hline
  6 &$\{(2,7),(3,7),(5,7)\}$ & $\frac{2}{3}(v_0-1)(v_0-2)^2(v_0-3)$ & $2(v_0-2)^2(v_0-3)$ &\\
  \hline
  7&$\{(1,3),(3,3),(3,5)\}$ & $5(v_0-1)(v_0-2)(v_0-3)$ & $10(v_0-2)(v_0-3)$ &\\
  \hline
  8&$\{(4,3),(6,3),(6,2)\}$ & $(v_0-1)(v_0-2)^2$ & $3(v_0-2)^2$ &\\
  \hline
  9&$\{(1,7),(2,1),(5,2)\}$ & $5(v_0-1)(v_0-2)^2$ & $10(v_0-2)^2$ &\\
  \hline
  10&$\{(1,4),(1,6),(4,5)\}$ & $\frac{1}{2}(v_0-1)(v_0-2)(v_0-3)(v_0-4) $ & $\frac{2}{3}(v_0-2)(v_0-3)(v_0-4) $ &\\
  \hline
  11&$\{(4,2),(4,6),(5,1)\}$ & $4(v_0-1)(v_0-2)^2$ & $12(v_0-2)^2$ &\\
  \hline
  12&$\{(2,7),(6,4),(6,7)\}$ & $3(v_0-1)(v_0-2)^2(v_0-3)$ & $14(v_0-2)^2(v_0-3)$ &\\
  \hline
  13&$\{(5,1),(6,1),(7,3)\}$ & $4(v_0-1)(v_0-2)^2(v_0-3)$ & $12(v_0-2)^2(v_0-3)$ &\\
  \hline
  14&$\{(2,2),(6,6),(7,1)\}$ & $3(v_0-1)(v_0-2)^2(v_0-3)$ & $14(v_0-2)^2(v_0-3)$ &\\
  \hline
  15&$\{(1,7),(5,4),(7,3)\}$ & $\frac{1}{2}(v_0-1)(v_0-2)^2(v_0-3)(v_0-4) $ & $\frac{7}{3}(v_0-2)^2(v_0-3)(v_0-4)$ &\\
  \hline
  16&$\{(3,5),(5,3),(5,4)\}$ & $\frac{1}{2}(v_0-1)(v_0-2)^2(v_0-3)(v_0-4) $ & $4(v_0-2)^2(v_0-3)(v_0-4)$ &\\
  \hline
  17&$\{(3,5),(5,7),(7,6)\}$ & $\frac{1}{6}(v_0-1)(v_0-2)^2(v_0-3)^2(v_0-4) $ & $3(v_0-2)^2(v_0-3)^2(v_0-4)$ &\\
  \hline
  18 &$\{(3,3),(5,6),(6,3)\}$ & $\frac{1}{2}(v_0-1)(v_0-2)^2(v_0-3)^2 $ & $4(v_0-2)^2(v_0-3)^2$ &\\
  \hline
  19&$\{(2,2),(3,5),(7,6)\}$ & $(v_0-1)(v_0-2)^2(v_0-3)^2$ & $8(v_0-2)^2(v_0-3)^2$ &\\
  \end{longtable}
Let $\mathcal{B}_i=C^G=\{(1,1),(1,2),a,b,c\}^G$, where $\{a,b,c\}$ is listed Case $i$ of Table \ref{BACaaa}. We claim that all blocks through points $(1,1)$ and $(1,2)$ have been considered. Evidently, for $1\le i<j\le 19$, we have   
 $ \mathcal{B}_i\cap \mathcal{B}_j=\emptyset$ and 
              \[ \sum_{\substack{1\le i\le 19\\\Phi \subset \mathcal{B}_i}} |\Phi|=\binom{v_0^2-2}{3}.\]
              
All pairs  $(|\Phi|,|\varPsi|)$ in Table \ref{BACaaa} can be calculated using the same method. To illustrate our approach clearly, we will use Case 7 of Table \ref{BACaaa} as an example for a detailed demonstration.
Let $\mathcal{B}_7=C^G$ where $C=\{(1,1),(1,2),(1,3),(3,3),(3,5)\}$.  
Every element of $\mathcal{B}_7$ has the form $$\{(h,i_1),(h,i_2),(h,i_3),(k,i_3),(k,i_4)\}, \{(i_1,h),(i_2,h),(i_3,h),(i_3,k),(i_4,k)\}.$$  
If $B\in \varPhi $ then $B$ can be classified into one of the following categories: 
\begin{flalign*}
 \{(1,1),(1,2),(1,m),(h,m),(h,n)\},&\quad\{(1,1),(1,2),(1,m),(h,1),(h,n)\},\\
 \{(1,1),(1,2),(1,m),(h,2),(h,n)\},&\quad\{(1,1),(1,2),(h,1),(h,m),(h,n)\},\\
\{(1,1),(1,2),(h,2),(h,m),(h,n)\},&\quad\{(1,1),(1,2),(h,1),(j,1),(k,2)\},\\
 \{(1,1),(1,2),(h,1),(j,2),(k,2)\}.&\quad                                  
\end{flalign*}
where $h, j, k\in \{2,3,4,\ldots,v_0\}$ (here $h, j, k$ are three distinct numbers) and $m, n\in \{3,4,\ldots,\\v_0\}$ with $m\ne n$.
Thus, $$|\varPhi |=3(v_0-1)(v_0-2)(v_0-3)+4(v_0-1)\binom{v_0-2}{2}=5(v_0-1)(v_0-2)(v_0-3).$$ 
If $B\in \varPsi $ then $B$ can be classified into one of the following categories:
\begin{flalign*}
  \{(1,1),(2,2),(1,2),(1,m),(2,n)\}, &\quad \{(1,1),(2,2),(1,m),(1,n),(2,n)\},  \\
  \{(1,1),(2,2),(1,m),(1,n),(2,n)\}, &\quad \{(1,1),(2,2),(1,m),(2,m),(2,n)\}, \\
  \{(1,1),(2,2),(2,1),(1,m),(1,n)\}, &\quad \{(1,1),(2,2),(1,2),(2,m),(2,n)\},  \\
  \{(1,1),(2,2),(2,1),(m,1),(n,2)\}, &\quad \{(1,1),(2,2),(m,1),(n,1),(n,2)\}, \\
  \{(1,1),(2,2),(1,2),(m,1),(n,2)\}, &\quad \{(1,1),(2,2),(m,1),(m,2),(n,2)\},\\
  \{(1,1),(2,2),(1,2),(m,1),(n,1)\}, &\quad \{(1,1),(2,2),(2,1),(m,2),(n,2)\}. 
\end{flalign*} 
where $m, n\in \{3,4,\ldots,v_0\}$ with $m\ne n$ .
Thus, $$|\varPsi |=8\binom{v_0-2}{1}\binom{v_0-3}{1}+4\binom{v_0-2}{2}=10(v_0-2)(v_0-3).$$

 Checking above 19 cases, we find that only cases 16-19 can deduce $|\varPhi|=|\varPsi|$ as $v_0\in\{7,9,19,39\}$. Accordingly, one of the following holds:
\begin{enumerate}
  \item[\rm(i)]  ${\rm Soc}(G)=A_9\times A_9$, and $\mathcal{D}$ is a $2$-$(81,5,\lambda)$ design where $\lambda \in\{5880,7056,14112\}$.
  \item[\rm(ii)]  ${\rm Soc}(G)=A_{19}\times A_{19}$, and $\mathcal{D}$ is a $2$-$(361,5,3329280)$ design.
\end{enumerate}

{\bf (b) $T$ is not 5-transitive group of $\Gamma$}

  Recall that $|\Gamma|=v_0\in\{7,9,19,39\}$ and $K$ is a primitive group (of almost simple or diagonal type) on $\Gamma$. Then $(v_0,T)=(7,PSL(2,7))$ or   $(9,PSL(2,8))$. Moreover,  if $T=PSL(2,7)$ then $G=PSL(2,7)\wr S_2$, and if $T=PSL(2,8)$ then $G$ is one of the following: $$PSL(2,8)\wr S_2,\, PSL(2,8)^2. S_3, \,PSL(2,8)^2. S_3, \, PSL(2,8)^2.6, \,P\varSigma L(2, 8)\wr S_2.$$
 Since $\cal D$ is a block-transitive $2$-$(v_0^2,5,\lambda)$, there must exist a base block $B$ of $\cal D$ that through points $(1,1)$ and $(1,2)$. Therefore, 
 \[B\in \bigcup_{\substack{1\le i\le 19\\\Phi \subset \mathcal{B}_i}}\Phi.\]
 By calculating the values of $|\Phi|$ and $|\Psi|$, we find that $G$ cannot be $PSL(2,7)\wr S_2$, furthermore, one of following holds:
\begin{enumerate} 
  \item[\rm(i)]  If $G=PSL(2,8)\wr S_2$, then $\mathcal{D}$ is a $2$-$(81,5,\lambda)$ design with $\lambda \in \{392,784,1568\}$.
  \item[\rm(ii)]  If $G=PSL(2,8)^2. S_3$ or $PSL(2,8)^2.6$, then $\mathcal{D}$ is a $2$-$(81,5,\lambda)$ design with $\lambda \in \{1176,2352,4704\}$.
  \item[\rm(iii)]  If $G=P\varSigma L(2, 8)\wr S_2$, then $\mathcal{D}$ is a $2$-$(81,5,\lambda)$ design with $\lambda \in \{1176,4704,7056\}$.
\end{enumerate}
Here, we have completed the proof of the Proposition \ref{pro5}. $\hfill\square$

Propositions \ref{pro3}-\ref{pro5} complete the proof of Theorem \ref{a1}. $\hfill\square$

\section*{Acknowledgements}
	This work is supported by the National Natural Science Foundation of China (No. 12361004) and the Natural Science Foundation of Jiangxi Province (Nos. 20224BAB211005 and 20242BAB25005). 

\section*{Conflict of interest statement}
The authors declare no conflict of interest.


\begin{thebibliography}{21}
\bibitem{Magma} W. Bosma, J. Cannon and C. Playoust: The  Magma algebra system I: The user language. J. Symb. Comput., 24 (1997), 235-265. 
\bibitem{Cameron} P. Cameron,  C. Praeger: Block-transitive designs I: point-imprimitive designs. Discrete Math., 118 (1993), 33-43. 
\bibitem{Handbook}  C. Colbourn, J. Dinitz: The CRC Handbook of Combinatorial Designs. CRC Press, Boca Raton, FL, 2007. 
\bibitem{Dav1} H. Davies:  Flag-transitivity and primitivity. Discrete Math., 63 (1987), 91-93. 
\bibitem{Delan} A. Delandtsheer, J. Doyen: Most block-transitive $t$-designs are point-primitive. Geom. Dedicata, 29 (1989), 307-310. 
\bibitem{Demb1968} P. Dembowski: Finite Geometries. Springer-Verlag, New York, 1968. 
\bibitem{Ding5} S. Ding, Y. Wang and X. Zhan: Flag-transitive $2$-$(v,5,\lambda)$ designs admitting a two-dimensional projective group. J. Math., 2024(1) (2024), 5521696. 
\bibitem{HMcL} D. Higman, J. McLaughlin: Geometric $ABA$-groups. Illinois J. Math., 5 (1961), 382-397. 
\bibitem{ONSC} M. Liebeck, C. Praeger and J. Saxl: On the O'Nan-Scott theorem for finite primitive permutation groups. J. Aust. Math. Soc. Ser. A, 44(3) (1988), 389-396. 
\bibitem{Montinaro} A. Montinaro: A classification of flag-transitive $2$-$(k^2,k,\lambda)$ designs with $\lambda|k$. J. Comb. Theory Ser. A, 197 (2023), 105750. 
\bibitem{Ray} D. Ray-Chaudhuri, R. Wilson, On $t$-designs. Osaka J. Math., 12 (1975), 737-744. 
\bibitem{Shen} J. Shen, S. Zhou: Flag-transitive $2$-$(v,5,\lambda)$ designs with sporadic socle. Front. Math. China, 15(6) (2020), 1201-1210. 
\bibitem{Zhu} D. Tian, S. Zhou: Flag-Transitive Point-Primitive Symmetric $(v,k,\lambda)$ Designs With $\lambda$ at Most $100$. J. Combin. Des., 21(4) (2013), 127-141. 
\bibitem{Tuan1} N. Tuan: A note on block-transitive point-imprimitive designs. European J. Combin., 24(1) (2003), 113-119. 
\bibitem{Tuan2} N. Tuan: On a divisibility problem for polynomials, and its application to Cameron-Praeger designs. J. Lond. Math. Soc., (2) 67(3) (2003), 545-560. 
\bibitem{Wie1964} H. Wielandt: Finite Permutation Groups. Academic Press, New York, 1964. 
\bibitem{ZhanDing}X. Zhan, S. Ding: A reduction for block-transitive triple systems. Discrete Math., 341 (2019), 2442-2447. 
\bibitem{ZhanDCC} X. Zhan, S. Ding and S. Bai: Flag-transitive 2-designs from $PSL(2,q)$ with block size 4. Des. Codes Cryptogr., 87(11) (2019), 2723-2728. 
\bibitem{ZhanPD} X. Zhan, X. Pang and S. Ding: Finite simple groups on triple systems. Ars Math. Contemp., 24(2) (2024), $\sharp$P2.09. 
\bibitem{ZhanZhou} X. Zhan, S. Zhou: A classification of flag-transitive point-imprimitive 2-designs with block size 6. J. Combin. Des., 26(4) (2018), 174-153. 
\bibitem{ZhanJCD}  X. Zhan, S. Zhou and G. Chen: Flag-transitive $2-(v,4,\lambda)$ designs of product type. J. Combin. Des., 26(9) (2018), 455-462. 
\bibitem{ZhanDM} X. Zhan, T. Zhou, S. Bai, S. Peng, L. Gan: Block-transitive automorphism groups on 2-designs with block size 4. Discrete Math., 343(7) (2020), 111726. 


\end{thebibliography}
\end{document}